\documentclass[runningheads]{llncs}

\usepackage{amssymb}
\usepackage{amsmath}
\usepackage{graphicx}
\usepackage{url}
\urldef{\mailsb}\path|jan.valdman@utia.cas.cz |

\providecommand{\keywords}[1]{\textbf{\textit{Keywords:}} #1}

\newcommand{\dx}{\,\mathrm{d}x}
\newcommand{\dy}{\,\mathrm{d}y}
\newcommand{\dxdy}{\dx \dy}

\begin{document}

\title{MATLAB Implementation of C1 finite elements: Bogner-Fox-Schmit rectangle\thanks{The work was supported by the Czech Science Foundation (GACR) through the grant GA18-03834S.} }
\titlerunning{MATLAB Implementation of C1 finite elements}
\author{Jan Valdman\inst{1,2}\orcidID{0000-0002-6081-5362} }
\authorrunning{Jan Valdman}
\institute{Institute of Mathematics, Faculty of Science, University of South Bohemia, 
Brani\v sovsk\' a 31, 37005~\v{C}esk\'{e}~Bud\v{e}jovice, Czech Republic \and
The Czech Academy of Sciences, 
Institute of Information Theory and Automation, 
Pod vod\'{a}renskou v\v{e}\v{z}\'{\i}~4, 18208~Praha~8, Czech Republic \\
\email{jan.valdman@utia.cas.cz}
}
\maketitle
\begin{abstract}
Rahman and Valdman (2013) introduced a new vectorized way to assemble finite element matrices. We utilize underlying vectorization concepts and extend MATLAB codes to implementation of Bogner-Fox-Schmit C1 rectangular elements in 2D. 
Our focus is on the detailed construction of elements and simple computer demonstrations including energies evaluations and their visualizations. \\ \\
\end{abstract}

\section{Introduction}
Boundary problems with fourth order elliptic operators 
\cite{Ciarlet-FEM}
appear in many applications including thin beams and plates and strain gradient elasticity \cite{Forest2009}. Weak formulations and implementations of these problems require H2-conforming finite elements, leading to $C^1$ continuity (of functions as well as of their gradients) of approximations over elements edges. This continuity condition is generally technical to achieve and few types of finite elements are known to quarantee it. We consider probably the simplest of them, the well known Bogner-Fox-Schmit rectangle \cite{BoFoSc65}, i.e., a rectangular $C^1$ element in two space dimensions.

We are primarily interested in explaining the construction of BFS elements, their practical visualization and evalutions.  Our MATLAB implementation is based on codes from \cite{AnjamValdman2015,RahmanValdman2013,HaVa14}. The main focus of these papers were assemblies of finite element matrices and local element matrices were computed all at once by array operations and stored in multi-dimentional arrays (matrices). Here, we utilize underlying vectorization concepts without the particular interest in corresponding FEM matrices.  More details on our recent implementations of $C^1$ models in nonlinear elastic models of solids can be found in \cite{KroemerValdman2019,FriedrichKruzikValdman2019,KroemerValdman2020}.
The complementary software to this paper is available at 
\begin{center}
\url{https://www.mathworks.com/matlabcentral/fileexchange/71346} 
\end{center}
for download and testing.

\section{Construction of $C^1$ finite elements}
\subsection{Hermite elements in 1D}
We define four cubic polynomials
\begin{equation}
\begin{split}
&\hat H_1(\hat x): = 2 \hat x^3 - 3 \hat x^2 + 1, \\
&\hat H_2(\hat x) : = -2 \hat x^3 + 3 \hat x^2, \\
&\hat H_3(\hat x) := \hat x^3 - 2 \hat x^2 + \hat x, \\
&\hat H_4(\hat x) := \hat x^3 - \hat x^2         
\end{split}
\end{equation}
over a reference interval $\hat I: = [0, 1]$
and can easily check the conditions:
\begin{equation}
\begin{split}
&\hat H_1(0) = 1, \quad \hat H_1(1) = 0, \quad \hat H_1'(0) = 0, \quad \hat H_1'(0) = 0, \\
&\hat H_2(0) = 0, \quad \hat H_2(1) = 1, \quad \hat H_2'(0) = 0, \quad \hat H_2'(0) = 0, \\
&\hat H_3(0) = 0, \quad \hat H_3(1) = 0, \quad \hat H_3'(0) = 1, \quad \hat H_3'(0) = 0, \\
&\hat H_4(0) = 0, \quad \hat H_4(1) = 0, \quad \hat H_4'(0) = 0, \quad \hat H_4'(0) = 1,
\end{split}
\label{conditions_1D}
\end{equation}
so only one value or derivative is equal to 1 and all other three values are equal to 0. These cubic functions create a finite element basis on $\hat I$. 
\begin{figure}[t]
  \centering
  \includegraphics[width=1\textwidth]{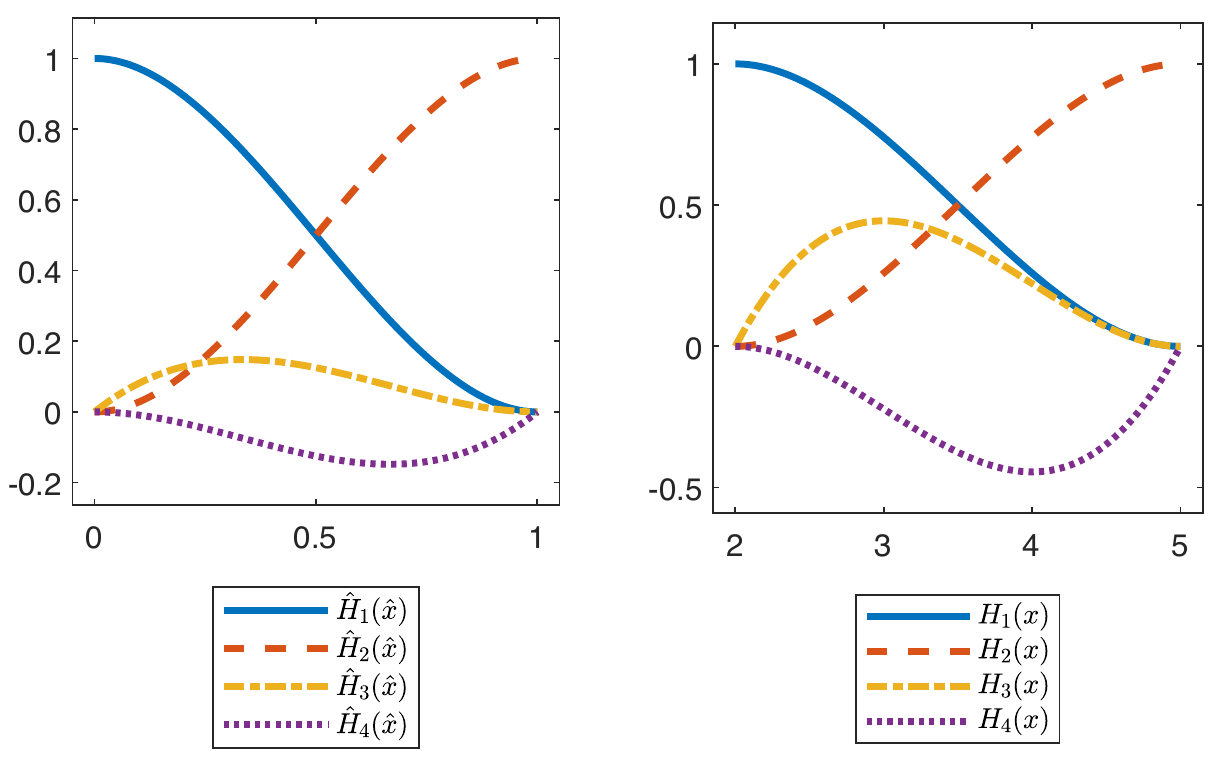}
  \caption{Reference basis functions $\hat H_i, i=1,\dots, 4$ on $[0,1]$ (left) and example of actual basis functions $H_i, i=1,\dots, 4$ on $[a, b]=[2, 5]$ (right).}
  \label{fig:C1_1D}
\end{figure}
More generally,
we define
\begin{equation}
\begin{split}
H_1(x) := \hat H_1(\hat x(x)), \\
H_2(x) := \hat H_2(\hat x(x)),  \\
H_3(\hat x) :=h \, \hat H_3(\hat x(x)),\\ H_4(\hat x) := h \, \hat H_4(\hat x(x))
\end{split}
\end{equation}
for $x \in I:=[a, b]$,
where $\hat x(x): =(x-a)/h$ is an affine mapping from $I$ to $\hat I$ and $h$ denotes the interval $I$ size 
$$h:=b-a.$$ These functions are also cubic polynomials and satisfy again the conditions \eqref{conditions_1D} with function arguments $0, 1$ replaced by $a, b$. They create actual finite element basis which ensures $C^1$ continuity of finite element approximations. The chain rule provides higher order derivatives:
\begin{equation}
\begin{split}
&H_1'(x)= \hat H_1'(\hat x) \, /h, \qquad H_1''(x)= \hat H_1''(\hat x) \, /h^2,\\
&H_2'(x)= \hat H_2'(\hat x) \, /h, \qquad H_2''(x)= \hat H_2''(\hat x) \, /h^2,\\
&H_3'(x)= \hat H_3'(\hat x) , \enspace \quad \qquad H_3''(x)= \hat H_3''(\hat x) \, /h,\\
&H_4'(x)= \hat H_4'(\hat x) ,  \enspace \quad \qquad H_4''(x)= \hat H_4''(\hat x) \, /h.\\
\end{split}
\label{1derivative_1D}
\end{equation}
\begin{example}
Example of basis functions defined on reference and actual intervals are shown in Figure \ref{fig:C1_1D} and pictures can be reproduced by 
\begin{verbatim}
  draw_C1basis_1D
\end{verbatim}
script located in the main folder. 
\end{example}

\subsection{Bogner-Fox-Schmit rectangular element in 2D}
\begin{figure}[t]
  \centering
  \includegraphics[width=1\textwidth]{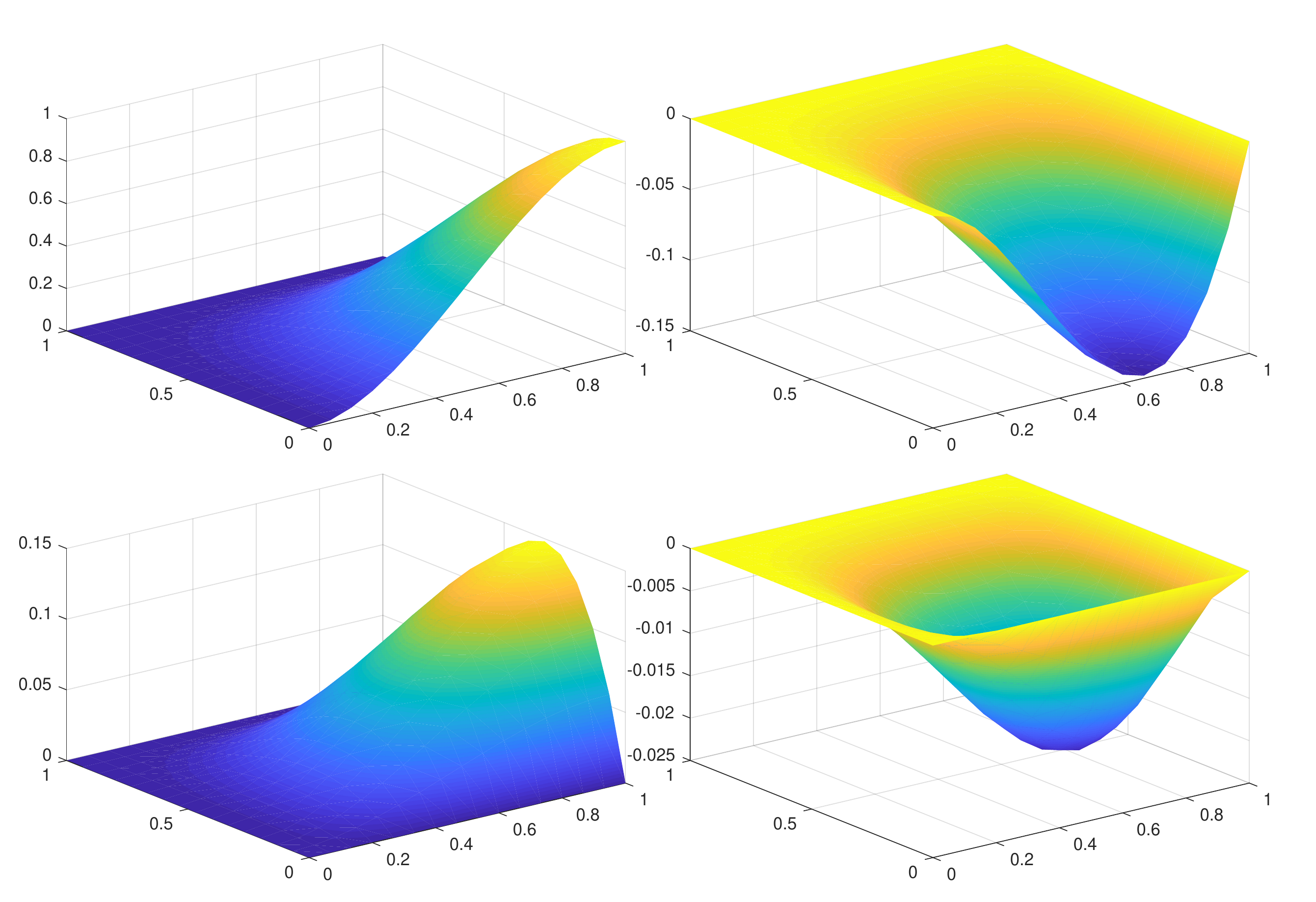}
  \caption{Bogner-Fox-Schmit basis functions 
  $\hat \varphi_{i}(\hat x, \hat y)$ for $i=2$ (top left), $i=6$ (top right), $i=8$ (bottom left), $i=13$ (bottom right) defined over a reference rectangle $\hat R =[0, 1] \times [0, 1]$.}
  \label{fig:C1_2D}
\end{figure}
Products of functions 
$$\tilde \varphi_{j,k}(\hat x, \hat y):=\hat H_j(\hat x) \, \hat H_k(\hat y), \quad j, k=1,\dots, 4$$
define 16 Bogner-Fox-Schmit (BFS) basis functions on a reference rectangle $\hat R :=[0 ,1] \times [0, 1]$.
For practical implementations, we reorder them as  
\begin{equation}
 \hat \varphi_{i}(\hat x, \hat y):= \tilde \varphi_{j_i,k_i}(\hat x, \hat y),
\quad i=1,\dots, 16,
\end{equation}
where sub-indices are ordered in a sequence 
\begin{equation}
\begin{split}
(j_i, k_i)_{i=1}^{16}
= \begin{matrix} \{
(1, 1), \, (2, 1),\, (2, 2),\, (1, 2), \\
\enspace (3, 1),\, (4, 1),\, (4, 2), \,(3, 2), \\
\enspace (1, 3), \,(2, 3), \,(2, 4), \,(1, 4), \\
\enspace \enspace(3, 3),\, (4, 3),\, (4, 4), \,(3, 4)
\}.
\end{matrix} 
\end{split}
\end{equation}
With this ordering, a finite element approximation $v \in C^1(\hat R)$ rewrites as a linear combination 
$$v(\hat x, \hat y) = \sum_{i=1}^{16} v_i \, \hat  \varphi_{i}(\hat x, \hat y),$$
where:
\begin{itemize}
    \item coefficients $v_1, \dots, v_4$ specify values of $v$,
    \item coefficients $v_5, \dots, v_8$ specify values of  $\frac{\partial v}{\partial x}$,
    \item coefficients $v_9, \dots, v_{12}$ specify values of   $\frac{\partial v}{\partial y}$,
    \item coefficients $v_{13}, \dots, v_{16}$ specify values of   $\frac{\partial^2 v}{\partial x \partial y}$ 
\end{itemize}
at nodes 
$$\hat N_1:=[0,0], \enspace \hat N_2:=[1,0], \enspace \hat N_3:=[1,1], \enspace \hat N_4:=[0,1].$$ 
\begin{example}
Four (out of 16) reference basis functions corresponding to the node $\hat N_2$ are shown in Figure \ref{fig:C1_2D} and pictures can be reproduced by 
\begin{verbatim}
  draw_C1basis_2D
\end{verbatim}
script located in the main folder.
\end{example}

\begin{figure}
  \centering
  \includegraphics[width=0.355\textwidth]{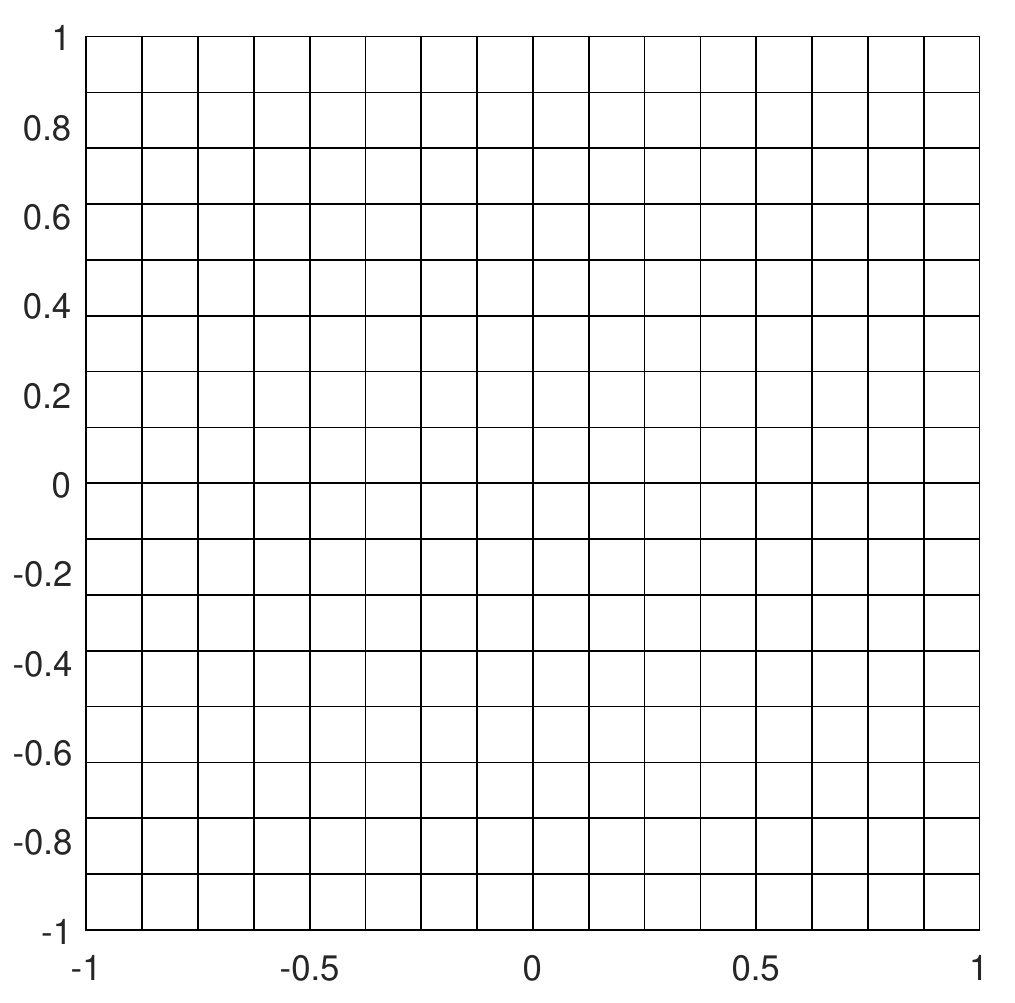}
  \hspace{0.06\textwidth}
  \includegraphics[width=0.56\textwidth]{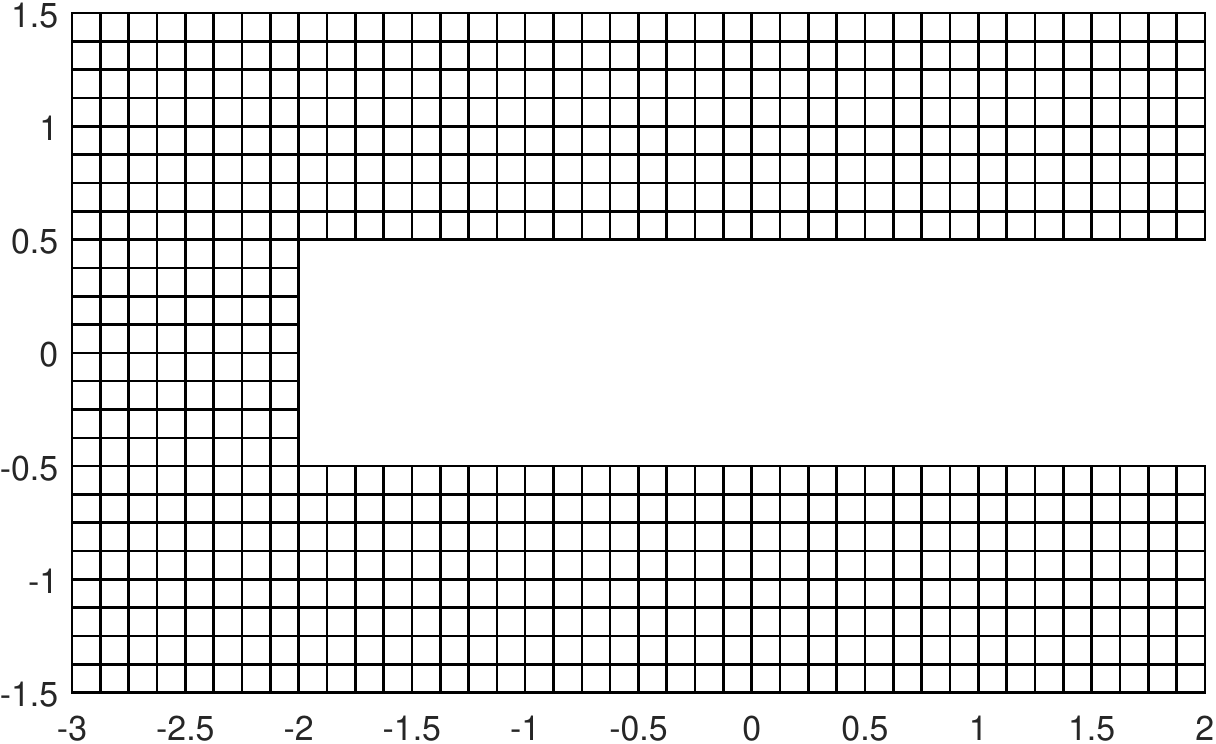}
  \caption{Examples of a triangulation in rectangles: of a square domain (left) and of a pincers domain (right) taken from \cite{KroemerValdman2019} and used for nonlinear elasticity simulations satisfying a non-selfpenetration condition. 
  }
  \label{fig:C1_2D_mesh}
\end{figure}

\subsubsection{More implementation details on functions evaluations}
For a general rectangle
$R:=[a, b] \times [c, d]$, we define an affine mapping
$$(\hat x, \hat y)(x, y) : =((x-a)/h_x, (y-c)/h_y),$$
from $R$ to $\hat R$, where the rectangular lengths are 
$$h_x:=b-a, \quad h_y:=d-c.$$ It enables us to define BFS basis functions on $R$ as
\begin{equation}
\varphi_{i}(x, y):= \hat H_{j_i} \left(\frac{x-a}{h_x} \right) \, \hat H_{k_i}\left(\frac{y-c}{h_y}\right),
\quad i=1,\dots, 16,
\label{BFS_basis}
\end{equation}
Based on \eqref{1derivative_1D}, higher order derivatives up to the second order, 
\begin{equation}\frac{\partial \varphi_{i}}{\partial x}, \enspace
   \frac{\partial \varphi_{i}}{\partial y}, \enspace
   \frac{\partial^2 \varphi_{i}}{\partial x^2}, \enspace
   \frac{\partial^2 \varphi_{i}}{\partial y^2}, \enspace   
   \frac{\partial^2 \varphi_{i}}{\partial x \partial y},
   \quad 
 i=1, \dots, 16
 \label{BFS_basis_derivatives}
 \end{equation}
can be derived as well. 
All basis functions \eqref{BFS_basis} are evaluated by the function
\begin{verbatim}
  shapefun(points',etype,h)
\end{verbatim}
and their derivatives  \eqref{BFS_basis_derivatives} by the function 
\begin{verbatim}
  shapeder(points',etype,h)
\end{verbatim}
For BFS elements, we have to set \verb+etype='C1'+ and a vector of rectangular lengths \verb+h=[hx,hy]+.
The matrix \verb+points+ then contains a set of points $\hat x \in \hat R$ in a reference element at which functions are evaluated. Both functions are vectorized 
and their outputs are stored as vectors, matrices or three-dimensional arrays.

 \begin{example}
The command 
\begin{verbatim}
  [shape]=shapefun([0.5 0.5]','C1',[1 1])
\end{verbatim}
 returns a (column) vector \verb+shape+ $\in \mathbb{R}^{16 \times 1}$ of all BFS basis function defined on $\hat R :=[0 ,1] \times [0, 1]$ and
 evaluated in the rectangular midpoint $ [0.5, 0.5] \in \hat R$. The command
\begin{verbatim}
  [dshape]=shapeder([0.5 0.5]','C1',[2 3])
\end{verbatim}
returns a three-dimensional array \verb+dshape+ $\in \mathbb{R}^{16 \times 1 \times 5}$ of all derivatives
 up to the second order of all BFS basis function defined on a general rectangle with lengths $h_x=2$ and $h_y=3$ and evaluated in the rectangular midpoint $ [0.5, 0.5] \in \hat R$.  
 
 For instance, if $R :=[1 ,3] \times [2, 5]$, values of all derivatives are evaluated in the rectangular midpoint $[2,3.5] \in R$.   
 
  More generally, if \verb+points+ $\in \mathbb{R}^{np \times 2}$ consists of $np>1$ points, then \verb+shapefun+ return a matrix of size $\mathbb{R}^{16 \times np}$ and 
\verb+shapeder+ returns a  three-dimensional array of size $\mathbb{R}^{16 \times np \times 5}$. 
\end{example}

\subsection{Representation of $C^1$ functions  and their visualization}

\begin{figure}[t]
  \centering
  \includegraphics[width=\textwidth,height=0.6\textheight]{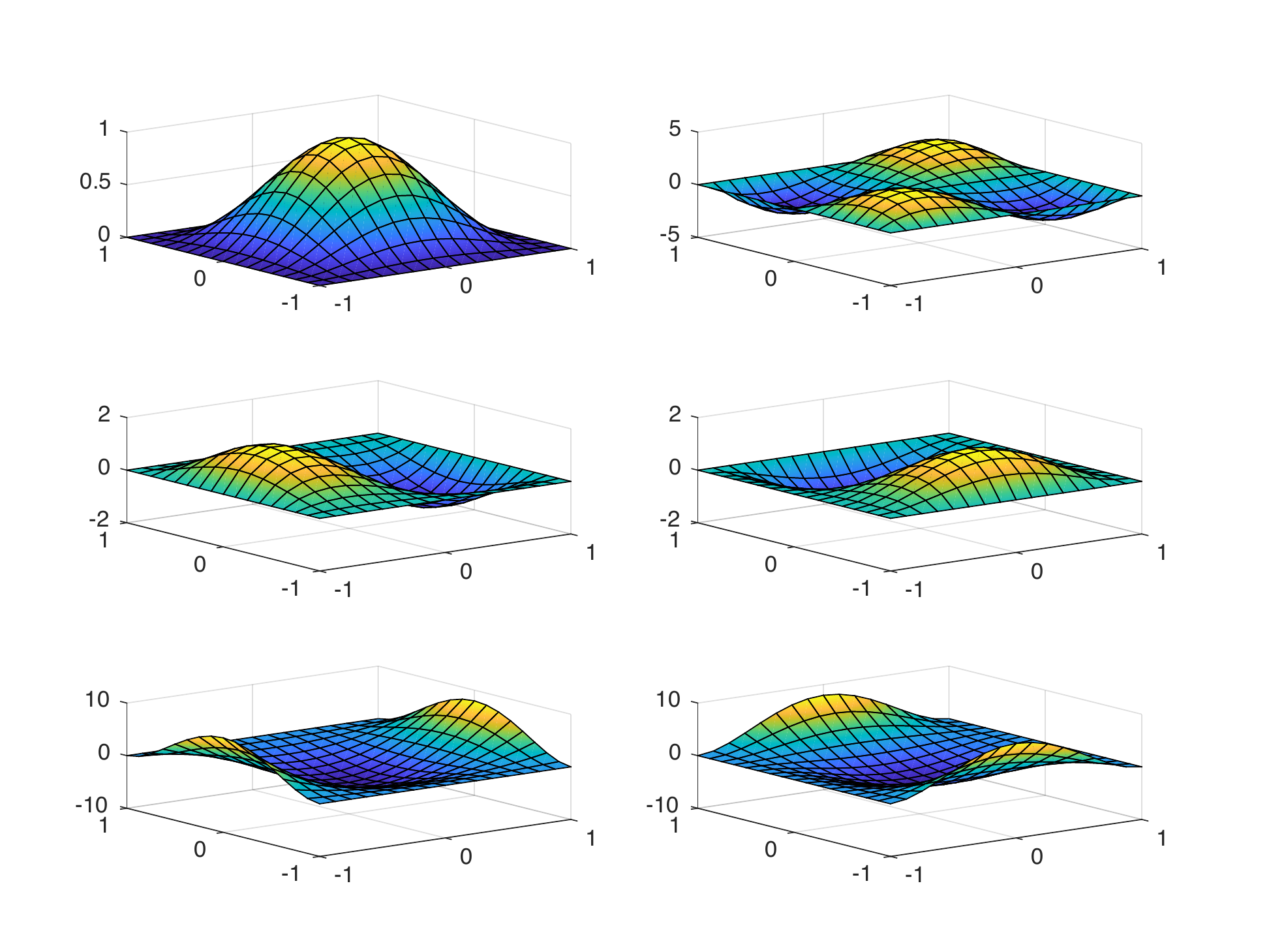}
  \caption{A function $v(x,y)=(1-x^2)^2(1-y^2)^2$ on $\Omega=(-1,1)^2$ represented in terms of BSF elements. Separate pictures show: 
  v (top left),
  $\frac{\partial^2 v}{\partial x \partial y}$ (top right),
  $\frac{\partial v}{\partial x}$ (middle left), 
  $\frac{\partial v}{\partial y}$ (middle right), 
  $\frac{\partial^2 v}{\partial x^2}$ (bottom left), $\frac{\partial^2 v}{\partial y^2}$ (bottom right). 
  }
  \label{fig:C1_2D_example}
\end{figure}
\begin{figure}[t]
  \centering
  \includegraphics[width=\textwidth]{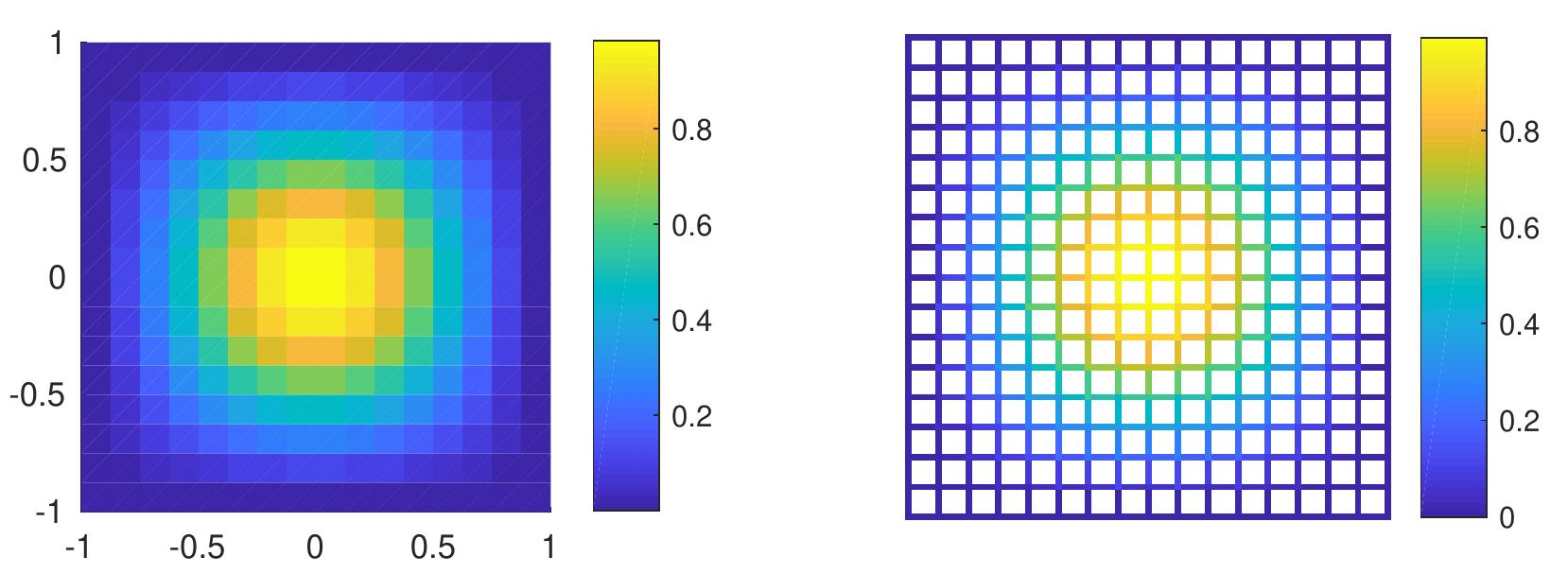}
  \caption{A function $v(x,y)=(1-x^2)^2(1-y^2)^2$ on $\Omega=(-1,1)^2$ and its values in elements midpoints (left) and edges midpoints (right). }
  \label{fig:densities}
\end{figure}

Let us assume a triangulation $\mathcal{T}(\Omega)$ into rectangles of a domain $\Omega$. In correspondence with our implementation, we additionally assume that all rectangles are of the same size, i.e., with lengths $h_x, h_y > 0$. 
Examples of $\mathcal{T}(\Omega)$ are given in Figure \ref{fig:C1_2D_mesh}.

Let $\mathcal{N}$ denotes the set of all rectangular nodes and $|\mathcal{N}|:=n$ the number of them. A $C^1$ function $v \in \mathcal{T}(\Omega)$ is represented in BSF basis by a matrix  
$$VC1=\begin{pmatrix} 
v(N_1), & \frac{\partial v}{\partial x}(N_1), & \frac{\partial v}{\partial y}(N_1),  & \frac{\partial^2 v}{\partial x \partial y}(N_1) \\
\vdots & \vdots & \vdots & \vdots \\
v(N_n), & \frac{\partial v}{\partial x}(N_n), & \frac{\partial v}{\partial y}(N_n),  & \frac{\partial^2 v}{\partial x \partial y}(N_n) \\
\end{pmatrix}$$
containing values of $v, \frac{\partial v}{\partial x}, \frac{\partial v}{\partial y}, \frac{\partial^2 v}{\partial x \partial y}$ in all nodes of $\mathcal{T}(\Omega)$. In the spirit of the finite element method, values of $v$ on each rectangle $T \in \mathcal{T}(\Omega)$ are obtained by an affine mapping to the reference element $\hat R$. Our implementations allows to evaluate and visualize continuous fields
$$ v, \, \frac{\partial v}{\partial x}, \, \frac{\partial v}{\partial y}, \,
\frac{\partial^2 v}{\partial x \partial y}$$
and also two additional (generally discontinuous) fields 
$$\frac{\partial^2 v}{\partial x^2}, \, \frac{\partial^2 v}{\partial y^2}.$$
It is easy to evaluate a $C^1$ function in a particular element point for all elements (rectangles) at once. A simple matrix-matrix MATLAB multiplication
\begin{verbatim}
  Vfun=VC1_elems*shape
\end{verbatim}
where a matrix \verb+VC1_elems+ $\in \mathbb{R}^{ne \times 16}$ contains in each row all 16 coefficients (taken from \verb+VC1+) corresponding to each element ($ne$ denotes a number of elements) returns a matrix \verb+Vfun+ $\in \mathbb{R}^{ne \times np}$
containing all function values in all elements and all points. Alternate multiplications
\begin{verbatim}  
   V1=VC1_elems*squeeze(dshape(:,1,:));                    % Dx
   V2=VC1_elems*squeeze(dshape(:,2,:));                    % Dy
  V11=VC1_elems*squeeze(dshape(:,3,:));                    % Dxx
  V22=VC1_elems*squeeze(dshape(:,4,:));                    % Dyy
  V12=VC1_elems*squeeze(dshape(:,5,:));                    % Dxy
\end{verbatim}
return matrices \verb+V1+, \verb+V2+, \verb+V11+, \verb+V22+, \verb+V12+ $\in \mathbb{R}^{ne \times np}$ containing values of 
all derivatives up to the second order in all elements and all points. A modification for evaluation of function values at particular edges points is also available and essential for instance for models with energies formulated on boundary edges \cite{KroemerValdman2020}.

\begin{example}
We consider a function 
\begin{equation}
 v(x,y)=(1-x^2)^2(1-y^2)^2   \label{v:example}
\end{equation}
on the domain $\Omega=(-1,1)^2$. This function was also used in \cite{FriedrichKruzikValdman2019} as an initial vertical displacement in a time-dependent simulation of viscous von K\'arm\'an plates. 

To represent $v$ in terms of BFS elements, we additionally need to know values of 
\begin{eqnarray}
\frac{\partial v}{\partial x}(x,y) = -4 x(1-x^2)(1-y^2)^2 \\
\frac{\partial v}{\partial x}(x,y) = -4 y(1-x^2)^2(1-y^2)\\
\frac{\partial^2 v}{\partial x \partial y}(x,y) = 16 x y(1-x^2)(1-y^2)
\end{eqnarray}
in nodes of a rectangular mesh $\mathcal{T}(\Omega)$. The function and its derivatives up to the second order are shown in Figure \ref{fig:C1_2D_example} and its values in elements and edges midpoints in Figure \ref{fig:densities}.

All pictures can be reproduced by 
\begin{verbatim}
  draw_C1example_2D
\end{verbatim}
script located in the main folder. 
\end{example}

\begin{figure}
  \centering
  \includegraphics[width=\textwidth]{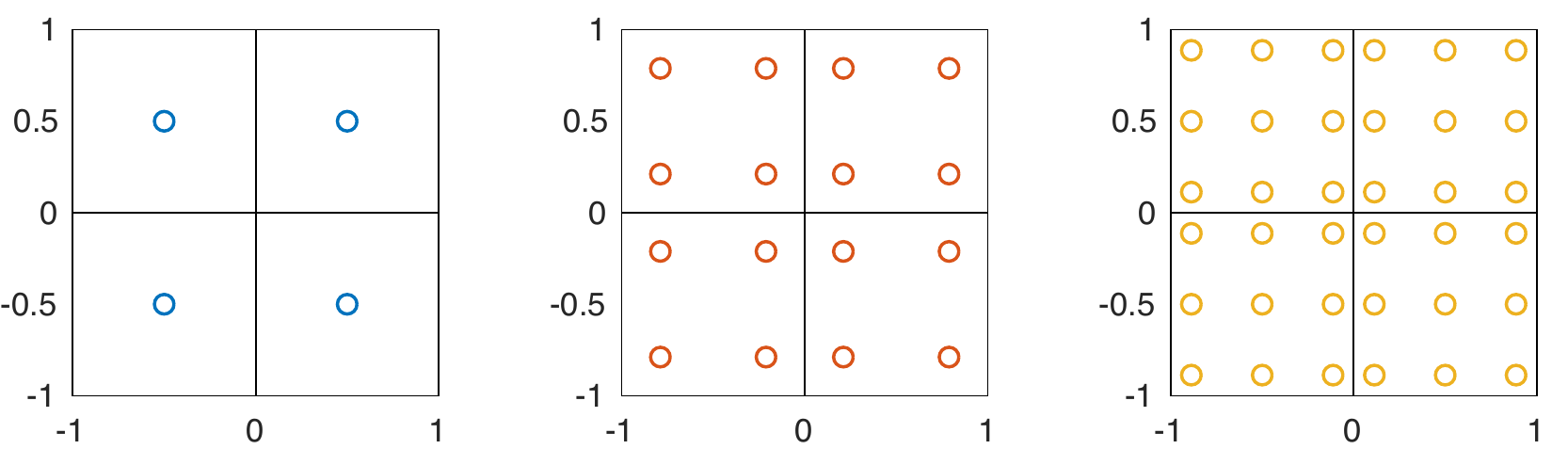}
  \caption{1, 4 and 9 Gauss points shown on actual rectangles of a square domain with 4 rectangles.}
  \label{fig:ips}
\end{figure}

\subsection{Evaluation of  $C^1$ functions and their numerical integration}
Various energy formulations include evaluations of integrals of the types
\begin{eqnarray}
 || v ||^2 :&=& \int_{\Omega} |v(x,y)|^2 \dxdy, \\
 || \nabla v ||^2 :&=& \int_{\Omega} |\nabla v(x,y)|^2 \dxdy, \\
  || \nabla^2 v ||^2 :&=& \int_{\Omega} |\nabla^2 v(x,y)|^2 \dxdy, \\
  (f, v) :&=& \int_{\Omega} f \, v \dxdy,
\end{eqnarray}
where $v \in H^2(\Omega)$ and $f \in L^2(\Omega)$ is given. The expression 
$$ (|| v ||^2 + || \nabla v ||^2 + || \nabla^2 v ||^2)^{1/2} $$ then defines the full norm in the Sobolev space $H^2(\Omega)$. For $v$ represented in the BFS basis we can evaluate above mentioned integrals numerically by quadrature rules. Our implementation provides three different rules with 1, 4 or 9 Gauss points. Each quadrature rule is defined by coordinates of all Gauss points and their weights.

\begin{figure}
  \centering
  \includegraphics[width=\textwidth]{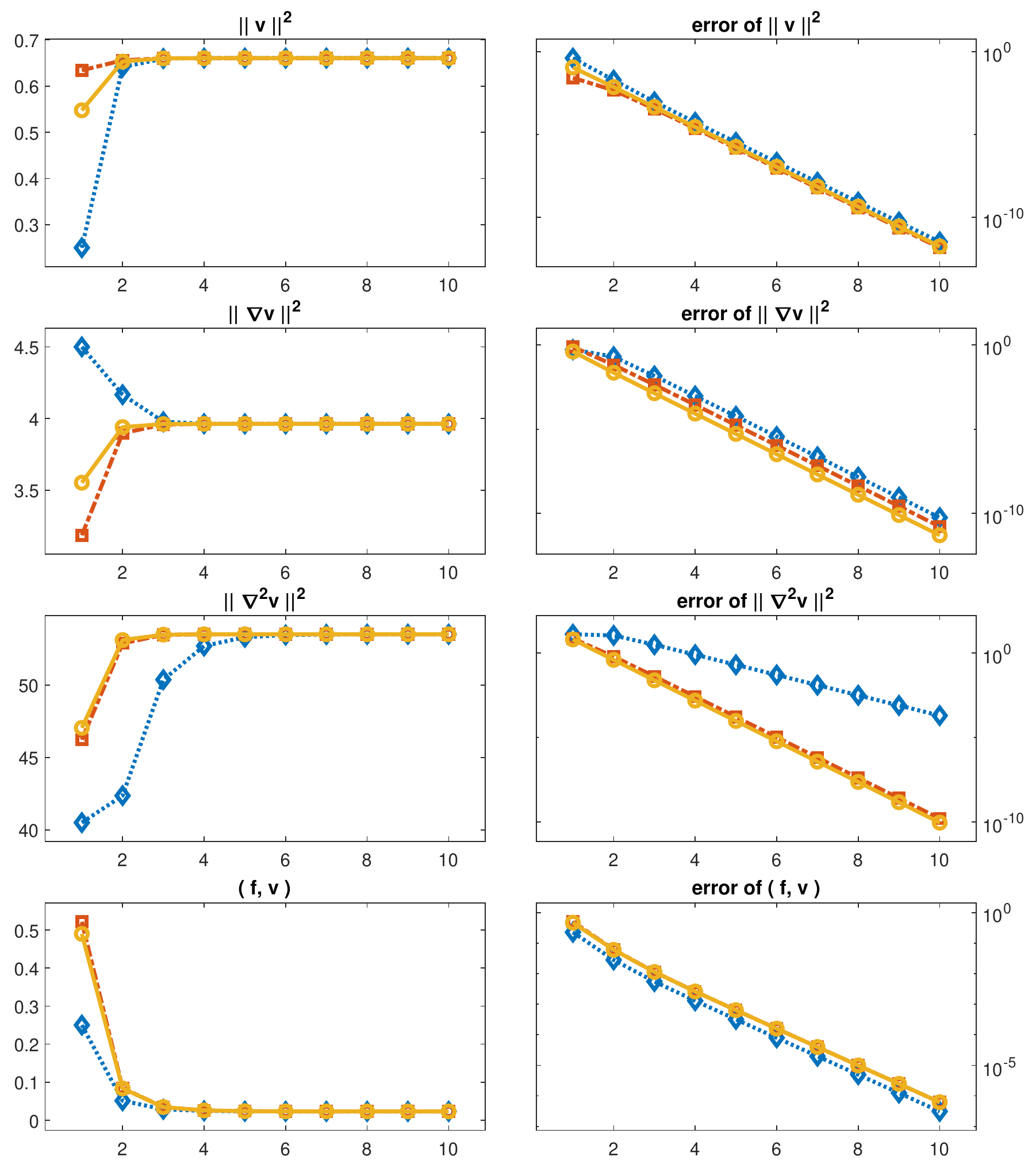}
  \caption{Values of integrals (the left column) and their absolute error (the right column) for levels 1-10 of uniform refinements using three different quadrature rules: 1 Gauss point - blue lines with diamonds, 4 Gauss points - yellow lines with circles, 9 Gauss points - red lines with squares.}
  \label{fig:convergence}
  \end{figure}
 
\begin{example}
Gauss points are displayed in Figure  \ref{fig:ips} and all pictures can be reproduced by 
\begin{verbatim}
  draw_ips
\end{verbatim}
script located in the main folder. 
\end{example}
\begin{example}
An analytical integration for the function $v$ from \eqref{v:example} and $f=x^2 y^2$
reveals that 
\begin{eqnarray*}
 || v ||^2 &=& 65536/99225 \approx 0.660478710002520, \\
 || \nabla v ||^2 &=& 131072/33075 \approx  3.962872260015117,\\
  || \nabla^2 v ||^2 &=& 65536/1225 \approx 53.498775510204084,\\
  (f, v) &=& 256/11025 \approx 0.023219954648526 
\end{eqnarray*}
for a domain $\Omega=(-1,1)^2$. We consider a sequence of uniformly refined meshes with levels 1-10:
\begin{itemize}
\item the coarsest (level=1) mesh with 9 nodes and 4 elements  is shown in Figure \ref{fig:ips},
\item a finer (level=4) mesh  with 289 nodes and 256 elements is shown in Figure \ref{fig:C1_2D_mesh} (left),
\item the finest (level=10) mesh consists of 1.050.625 nodes and 1.048.576 elements, not shown here.
\end{itemize}
Figure \ref{fig:convergence}  the depicts the convergence of numerical quadratures to the exact values above. 
 We notice that all three quadrature rules provide the same rates of convergence. The only exception is the evaluation of the second gradient integral $|| \nabla^2 v ||^2$, where the numerical quadrature using 1 Gauss point  deteriorates the convergence. 
 
All pictures can be reproduced by 
\begin{verbatim}
  start_integrate
\end{verbatim}
script located in the main folder. 
\end{example}

%
%
%
%

\bibliographystyle{splncs03}

\end{document}